\documentstyle[amsfonts,12pt]{article}
\newcommand{\f}{{\cal O}}
\newcommand{\cD}{{\cal D}}
\renewcommand{\P}{{\Bbb P}}
\newcommand{\C}{{\Bbb C}}
\newcommand{\R}{{\Bbb R}}

\newcommand{\om}{\omega}

\newcommand{\g}{{\frak g}}

\newcommand{\e}{{\frak e}}
\newcommand{\gl}{{\frak gl}}

\newcommand{\SL}{\mbox{SL}}
\newcommand{\SO}{\mbox{SO}}

\newcommand{\E}{\mbox{E}}

\newcommand{\T}{\mbox{T}}
\renewcommand{\H}{\mbox{H}}
\newcommand{\eH}{\mbox{\em H}}

\newcommand{\Ad}{\mbox{Ad}}

\newcommand{\rank}{\mbox{rank\ }}
\newcommand{\lon}{\mbox{$\longrightarrow$}}
\newcommand{\rar}{\rightarrow}
\newcommand{\hook}{\hookrightarrow}
\newcommand{\p}{\partial}
\newcommand{\ot}{\otimes}

\newtheorem{theorem}{Theorem}

\newtheorem{prop}{Proposition}

\newtheorem{lemma}[theorem]{Lemma}

\newcommand{\contr}{\,\begin{picture}(5,5)
\put(0,0){\line(5,0){5}}
\put(5,0){\line(0,0){5}} 
\end{picture}\,}     
\newcommand{\TXX}{\begin{picture}(75,20)
                 \put(0,2){\line(5,0){75}}
                 \put(30,-17){\line(0,0){17}}
                 \put(0,2){\circle*{4}} 
                 \put(15,2){\circle*{4}}                 
                 \put(30,2){\circle*{4}} 
                 \put(45,2){\circle*{4}} 
                 \put(60,2){\circle*{4}} 
                 \put(71,-1){$\times$} 
                 \put(30,-17){\circle*{4}} 
                 \put(-2,6){\shortstack{1}}
                 \put(13,6){\shortstack{0}}
                 \put(28,6){\shortstack{0}}
                 \put(43,6){\shortstack{0}}
                 \put(58,6){\shortstack{0}}
                 \put(73,6){\shortstack{0}}
                 \put(34,-21){\shortstack{0}}  
                 \end{picture}}
\newcommand{\TT}{\begin{picture}(75,20)
                 \put(0,2){\line(5,0){75}}
                 \put(30,-17){\line(0,0){17}}
                 \put(0,2){\circle*{4}} 
                 \put(15,2){\circle*{4}}                 
                 \put(30,2){\circle*{4}} 
                 \put(45,2){\circle*{4}} 
                 \put(60,2){\circle*{4}} 
                 \put(75,2) {\circle*{4}} 
                 \put(30,-17){\circle*{4}} 
                 \put(-2,6){\shortstack{1}}
                 \put(13,6){\shortstack{0}}
                 \put(28,6){\shortstack{0}}
                 \put(43,6){\shortstack{0}}
                 \put(58,6){\shortstack{0}}
                 \put(73,6){\shortstack{0}}
                 \put(34,-21){\shortstack{0}}  
                 \end{picture}}
\newcommand{\LL}{\begin{picture}(75,20)
                 \put(0,2){\line(5,0){75}}
                 \put(30,-17){\line(0,0){17}}
                 \put(0,2){\circle*{4}} 
                 \put(15,2){\circle*{4}}                 
                 \put(30,2){\circle*{4}} 
                 \put(45,2){\circle*{4}} 
                 \put(60,2){\circle*{4}} 
                 \put(71,-1){$\times$} 
                 \put(30,-17){\circle*{4}} 
                 \put(-2,6){\shortstack{0}}
                 \put(13,6){\shortstack{0}}
                 \put(28,6){\shortstack{0}}
                 \put(43,6){\shortstack{0}}
                 \put(58,6){\shortstack{0}}
                 \put(73,6){\shortstack{1}}
                 \put(34,-21){\shortstack{0}}  
                 \end{picture}}
\newcommand{\TXO}{\begin{picture}(75,20)
                 \put(0,2){\line(5,0){75}}
                 \put(30,-17){\line(0,0){17}}
                 \put(0,2){\circle*{4}} 
                 \put(15,2){\circle*{4}}                 
                 \put(30,2){\circle*{4}} 
                 \put(45,2){\circle*{4}} 
                 \put(60,2){\circle*{4}} 
                 \put(71,-1){$\times$} 
                 \put(30,-17){\circle*{4}} 
                 \put(-2,6){\shortstack{2}}
                 \put(13,6){\shortstack{0}}
                 \put(28,6){\shortstack{0}}
                 \put(43,6){\shortstack{0}}
                 \put(58,6){\shortstack{0}}
                 \put(72,6){\shortstack{-1}}
                 \put(34,-21){\shortstack{0}}  
                 \end{picture}}
\newcommand{\TXP}{\begin{picture}(75,20)
                 \put(0,2){\line(5,0){75}}
                 \put(30,-17){\line(0,0){17}}
                 \put(0,2){\circle*{4}} 
                 \put(15,2){\circle*{4}}                 
                 \put(30,2){\circle*{4}} 
                 \put(45,2){\circle*{4}} 
                 \put(60,2){\circle*{4}} 
                 \put(71,-1){$\times$} 
                 \put(30,-17){\circle*{4}} 
                 \put(-2,6){\shortstack{0}}
                 \put(13,6){\shortstack{0}}
                 \put(28,6){\shortstack{0}}
                 \put(43,6){\shortstack{0}}
                 \put(58,6){\shortstack{1}}
                 \put(71,6){\shortstack{-1}}
                 \put(34,-21){\shortstack{0}}  
                 \end{picture}}
\newcommand{\TXL}{\begin{picture}(75,20)
                 \put(0,2){\line(5,0){75}}
                 \put(30,-17){\line(0,0){17}}
                 \put(0,2){\circle*{4}} 
                 \put(15,2){\circle*{4}}                 
                 \put(30,2){\circle*{4}} 
                 \put(45,2){\circle*{4}} 
                 \put(60,2){\circle*{4}} 
                 \put(71,-1){$\times$} 
                 \put(30,-17){\circle*{4}} 
                 \put(-2,6){\shortstack{3}}
                 \put(13,6){\shortstack{0}}
                 \put(28,6){\shortstack{0}}
                 \put(43,6){\shortstack{0}}
                 \put(58,6){\shortstack{0}}
                 \put(71,6){\shortstack{-2}}
                 \put(34,-21){\shortstack{0}}  
                 \end{picture}}
\newcommand{\TXM}{\begin{picture}(75,20)
                 \put(0,2){\line(5,0){75}}
                 \put(30,-17){\line(0,0){17}}
                 \put(0,2){\circle*{4}} 
                 \put(15,2){\circle*{4}}                 
                 \put(30,2){\circle*{4}} 
                 \put(45,2){\circle*{4}} 
                 \put(60,2){\circle*{4}} 
                 \put(71,-1){$\times$} 
                 \put(30,-17){\circle*{4}} 
                 \put(-2,6){\shortstack{1}}
                 \put(13,6){\shortstack{0}}
                 \put(28,6){\shortstack{0}}
                 \put(43,6){\shortstack{0}}
                 \put(58,6){\shortstack{1}}
                 \put(71,6){\shortstack{-2}}
                 \put(34,-21){\shortstack{0}}  
                 \end{picture}}
\newcommand{\TXN}{\begin{picture}(75,20)
                 \put(0,2){\line(5,0){75}}
                 \put(30,-17){\line(0,0){17}}
                 \put(0,2){\circle*{4}} 
                 \put(15,2){\circle*{4}}                 
                 \put(30,2){\circle*{4}} 
                 \put(45,2){\circle*{4}} 
                 \put(60,2){\circle*{4}} 
                 \put(71,-1){$\times$} 
                 \put(30,-17){\circle*{4}} 
                 \put(-2,6){\shortstack{0}}
                 \put(13,6){\shortstack{0}}
                 \put(28,6){\shortstack{0}}
                 \put(43,6){\shortstack{0}}
                 \put(58,6){\shortstack{0}}
                 \put(72,6){\shortstack{0}}
                 \put(34,-21){\shortstack{0}}  
                 \end{picture}}

\begin{document}

\title{Exotic holonomies $\E_7^{(a)}$}

\author{Q.-S.\ Chi, S.A.\ Merkulov, L.J.\ Schwachh\"ofer}

\date{ }
\maketitle

\begin{abstract}
It is proved that the Lie groups  
$\E_7^{(5)}$ and $\E^{(7)}_7$ represented in $\R^{56}$ and 
the Lie group $\E_7^{\C}$ represented in $\R^{112}$
occur as holonomies of torsion-free affine connections.
It is also shown that the moduli spaces of torsion-free affine
connnections with these holonomies are finite dimensional, and
that every such connection has a local symmetry group of positive
dimension.
\end{abstract}


\begin{center}
{\bf \S 1  Introduction}
\end{center}
The notion of the {\em Holonomy} of an affine connection was
introduced originally by \'Elie Cartan in the 1920s who used it as 
an important tool in his attempt to classify all locally 
symmetric manifolds. Over time, the holonomy group proved to be one
of the 
most informative and useful characteristics of an affine connection
and found 
many applications in both mathematics and physics.

By definition, the holonomy of an affine connection on a connected
manifold $M$ is the subgroup of all linear automorphisms
of $T_pM$ which are induced by parallel translation along $p$-based
loops.

In 1955, Berger \cite{Berger} showed that the list of irreducibly 
acting matrix Lie groups which can, in principle, occur as the 
holonomy of a torsion-free affine connection is very restricted. 
Berger presented his classification of all possible candidates for
irreducible holonomies in two parts. The first part contains
all possible groups preserving a non-degenerate 
symmetric bilinear form, the second part consists of those groups
which do {\em not} preserve such a form; the latter part was stated to
be complete up to a finite number of missing terms and was given without
a proof.

Bryant \cite{Bryant2} was the first to discover the incompleteness of
the second part of Berger's list, and referred to the missing entries as
{\em exotic 
holonomies}. Since then, several other families of exotic holonomies
have been found
\cite{Bryant3,CS,CMS}. In this paper we present one more family of 
exotic holonomies associated with various real forms of the
complex 56-dimensional representation of $\E_7^{\C}$.

\paragraph{Main Theorem.}{\em
\begin{itemize}
\item[(i)] All representations in the following
table

{\em
\begin{center}
\begin{tabular}{|c|c|c|c|}
\hline
\begin{tabular}{c} \vspace{-3mm} \\ Group: \\ \vspace{-3 mm} \end{tabular}
 & $\E_7^{(5)}$ &  $\E_7^{(7)}$ & $\E_7^{\C}$ \\
\hline
\begin{tabular}{c}\vspace{-3mm}\\ Representation space: \\ 
\vspace{-3mm} \end{tabular}
 & $\R^{56}$ & $\R^{56}$ & $\R^{112}$ \\
\hline
\end{tabular}
\end{center}
}

occur as holonomies of torsion-free affine connections.

\item[(ii)] Any torsion-free affine connection with one of these holonomies
is analytic.

\item[(iii)] The moduli space of torsion-free affine connections 
with one of these holonomies is finite dimensional.

\item[(iv)] Any such connection has a (local) symmetry group of positive
dimension.
\end{itemize}
}

This theorem is proved by combining twistor techniques of 
\cite{Me1} used to compute all the necessary $\E_7$-modules
$K(\e_7)$, $K^1(\e_7)$ and ${\cal P}^{(1)}(\e_7)$ with the
construction of torsion-free affine connections with prescribed holonomy
via deformations of a certain linear Poisson structure \cite{CMS}.

\begin{center}
{\bf \S 2 \, Borel-Weil approach to $\E_7^{(a)}$ } 
\end{center}
Let $V$ be a vector space and $\g$ an irreducible Lie subalgebra of 
$\gl(V)\simeq V\ot V^*$. In the holonomy group  context, 
one is interested in the following three $\g$-modules:
\begin{itemize}
\item[(i)] $\g^{(1)}:= (\g\ot V^*)\cap (V\ot \odot^2 V^*$),  
\item[(ii)] the {\em curvature space} $K(\g):=\ker i_1$, where $i_1$
is the composition
$$
i_1: \g\ot \Lambda^2 V^* \lon V\ot V^*\ot \Lambda^2 V^* \lon
V\ot \Lambda^3 V^* ,
$$
\item[(iii)] the {\em 2nd curvature space} $K^1(\g):=\ker i_2$, where
$i_2$ is the composition 
$$
i_2: K(\g)\ot V^* \lon \g\ot \Lambda^2 V^* \ot V^* \lon V \ot\Lambda^3 V^*.
$$
\end{itemize}
Note that if $\p$ is the composition
$$
\p: \g^{(1)}\ot V^* \rar \g\ot V^*\ot V^* \rar \g\ot \Lambda^2 V^*
$$
then $\p(\g^{(1)}\ot V^*)\subset K(\g)$.

The geometric meaning of $\g^{(1)}$ is that if there exists a (local) 
torsion-free affine connection $\nabla$ on a manifold $M$ with holonomy 
algebra $\g$ then, for any (local) function $\Gamma: M\rar \g^{(1)}$, the
affine connection $\nabla + \Gamma$ is again torsion-free and has holonomy
algebra $\g$; thus, in some sense, $\g^{(1)}$ measures the non-uniqueness
of torsion-free affine connections with holonomy $\g$ on a fixed
manifold.

The significance of $K(\g)$ and $K^1(\g)$ is that 
the curvature tensor (the covariant derivative of the curvature tensor
respectively) of a torsion-free affine connection
$\nabla$ with holonomy $\g$ at a point $p \in M$ is represented by an
element of $K(\g)$ ($K^1(\g)$ respectively).

Therefore, $\g$ can be a candidate to the holonomy algebra of a 
torsion-free affine connection only if $K(\g)\neq 0$. The question 
then remains how to compute $K(\g)$.

With any real irreducible representation of a real reductive Lie algebra
one may associate an irreducible complex representation of a complex
reductive Lie algebra. Since all the above $\g$-modules behave reasonably
well under this association, we may assume from now on that
$V$ is a finite dimensional complex vector space and $\g\subset
\gl(V)$ is an irreducible representation of a complex reductive Lie
algebra. Clearly, $G=\exp(\g)$ acts irreducibly in $V^*$ via the dual
representation. Let $\tilde{X}$ be the  $G$-orbit of a highest weight 
vector in $V^*\setminus 0$. Then the quotient $X:= \tilde{X}/{\Bbb C}^*$
is a  compact complex homogeneous-rational manifold  canonically  
embedded into ${\P}(V^*)$, and there is a commutative diagram
$$
\begin{array}{ccc}
\tilde{X} & \hook & V^* \setminus 0 \\
   \downarrow &       & \downarrow \\
        X   & \hook &   \P(V^*) 
\end{array}
$$
In fact, $X=G_s/ P $, where $G_s$ is the semisimple part of $G$ and $P$
is the parabolic subgroup of $G_s$ leaving  a highest
weight vector in $V^*$ invariant up to a scalar. Let $L$ be the restriction of
the hyperplane section bundle $\f(1)$ on $\P(V^*)$ to the submanifold $X$.
Clearly, $L$ is an ample  homogeneous line bundle on $X$.
We call $(X,L)$ the {\em Borel-Weil data}\, associated with $(\g,V)$.

According to Borel-Weil, the representation space
$V$ can be easily reconstructed from $(X,L)$ as $V = \H^0(X,L)$.
What about $\g$? The Lie algebra of
the Lie group of all global biholomorphisms of the line bundle $L$
which commute with the projection $L\rar X$ is isomorphic
to $H^0(X,L\ot (J^1 L)^*)$ --- a central extension of the Lie
algebra $\H^0(X,TX)$. Whence, as a complex Lie algebra, 
$\H^0(X,L\ot (J^1 L)^*)$ has a natural complex irreducible 
representation in $\H^0(X,L)=V$; with very few (and well studied in the 
holonomy context) exceptions \cite{Ah}, this representation is,
up to a central extension, isomorphic to the original $\g$.

Remarkably enough, the basic $\g$-modules  defined above fit nicely into
the Borel-Weil pa\-ra\-digm as well.

\begin{prop}{\em\cite{Me1}}\label{Spencer}
For a compact complex homogeneous-rational manifold $X$ and an
ample line bundle $L\rar X$, there is an isomorphism
$$
\g^{(1)} = \eH^0\left(X, L\otimes \odot^{2} N^*\right), 
$$ 
and  an exact sequence of $\g$-modules,
$$
0\lon \frac{K(\g)}{\p(\g^{(1)}\ot V^*)}  
\lon \eH^1\left(X,L\ot \odot^{3} N^*\right)
\lon \eH^1\left(X,L\ot
\odot^{2}N^*\right) \ot V^*   ,                
$$
where $\g$ is $\eH^0(X,L\otimes N^*)$ represented in $V=\eH^0(X,L)$.
\end{prop}
{\em Proof}. The result follows easily from the 
exact sequences
$$
0\lon L\ot \odot^2 N^* \lon L\ot N^* \ot V^* \lon L\ot \Lambda^2V^* 
$$
and
$$
0\lon L\ot \odot^3 N^* \lon L\ot \odot^2 N^* \ot V^* \lon L\ot
N^*\ot \Lambda^2V^* \lon L\ot \Lambda^3V^*,  
$$
where arrows are a combination of a 
natural monomorphism $N^* \lon V^*\ot \f_X$
(which holds due to ampleness of $L$)  with the
antisymmetrization. $\Box$

It is well known that the complex exceptional Lie algebra $\e_7^{\C}$
has four real forms $\e_7^{(4)}$, $\e_7^{(5)}$, $\e_7^{(6)}$ and
$\e_7^{(7)}$ with signatures 0, 54, 64 and 70 respectively 
(see, e.g., \cite{HS,On}). Two of these, $\e_7^{(5)}$ and $\e_7^{(6)}$,
can be irreducibly represented in $\R^{56}$.
 Let $\rho$ denote  the irreducible real 
representation $\e_7^{(a)}\rar\gl(V)$, where $V$ is $\R^{56}$ for 
$a=5,7$ and $V=\C^{56}\simeq 
\R^{112}$ for $(a)=\C$. Let $\Ad: \e_7^{(a)}\rar \gl(\e_7^{(a)})$ 
denote the adjoint representation.

\begin{theorem}\label{T1}
$$
K(\rho(\e_7^{(a)})) \simeq {\em\Ad} (\e_7^{(a)}), \ \ \ \ \
K^1(\rho(\e_7^{(a)})) \simeq V^*.
$$
\end{theorem}  
{\em Proof}. We shall prove this statement for the complex
representation only. That it is true for real representations as well
will follow from the invariance of all the constructions under
the associated real structures in $e_7^{\C}$. 

 Let $(X,L)$ be the Borel-Weil data associated to $\rho:\e_7^{\C}
\rar \gl(V)$, $V\simeq\C^{56}$. Then $X=\E_7^{\C}/P$ is a 27-dimensional
compact 
complex homogeneous-rational manifold whose tangent bundle has, as an
irreducible homogeneous vector bundle, the Dynkin diagram
representation \cite{BE}
$$
TX=\ \ \TXX \ \ \ ,   \vspace{7 mm}
$$
while $L\rar X$ is given by
$$
\ \ \, L=\ \ \LL \ \ \ .  \vspace{7 mm}
$$
Here and below the weights of irreducible homogeneous vector bundles 
are given in the basis of fundamental weights. 

Using Konstant's formula and Table 5 in the reference chapter of
 \cite{On} to find irreducible decompositions of tensor powers of the 
simplest 27-dimensional irreducible representation of $E_{6}^{\C}$
(which, in our case, is isomorphic to the semisimple part of the parabolic
$P$), one obtains
$$
\odot^2 TX\ot L^* = \ \ \TXO \ \ \ + \ \ \ \TXP \vspace{7 mm}
$$
and  
$$
\odot^3 TX\ot L^{*2} = \ \ \, \TXL \ \ \, + \ \ \, \TXM \ \ \, + \ \,
 \ \TXN.  \vspace{7 mm}
$$
Then, using the long exact sequences of the extensions
$$
0\lon N^* \lon L\ot \odot^2 N^* \lon \odot^2 TX \ot L^* \lon 0,
$$
$$
0\lon \odot^2 N^* \lon L\ot \odot^3 N^* \lon \odot^3 TX \ot L^{2*} \lon 0,
$$
the Bott-Borel-Weil theorem and Proposition~\ref{Spencer},
one easily finds
$$
\H^0(X,L\ot \odot^2 N^*) = \H^1(X,L\ot \odot^2 N^*) = 0
$$
and
$$
K(\rho(\e_7^{\C} + \C)) \simeq \H^1(X,L\ot \odot^3 N^*)=
\ \ \TT \ \ \simeq \Ad(\e_{7}^{\C}). \vspace{7 mm}
$$

Let us find next the explicit form of $K(\rho(\e_7^{\C} + \C))$
as a subset of all elements  in $\rho(\e_7^{\C} + \C))\ot \Lambda^2 V^*$
satisfying the first Bianchi identities.

Recall \cite{Adams} that $\rho: \e_7^{\C} \rar \gl(V)$ enjoys a
non-zero invariant skew symmetric invariant product
\begin{eqnarray*}
\Lambda^2 V & \lon & \C \\
    u\ot v  & \lon & <\,u,v\,>,
\end{eqnarray*}
and a non-zero invariant symmetric map
\begin{eqnarray*}
\odot^2 V & \lon & e_7^{\C} \\
    u\ot v  & \lon & u\circ v,
\end{eqnarray*}
which are unique up to non-zero scalar factor and satisfy
\begin{eqnarray*}
<\rho(A)u, v> &=& \lambda B(A,u\circ v)\\
B(u\circ v, s\circ t) - B(u\circ t, s\circ v)&=& \mu
\left(2<\!u,s\!><\!v,t\!> - <\!u,t\!><\!v,s\!> - <\!u,v\!><\!s,t\!>\right)
\end{eqnarray*}
for all $A\in e_7^{\C}$ and $u,v,s,t\in V$. Here $\lambda$ and $\mu$
are fixed non-zero constants and $B(\ , \ )$ is the Killing form.

Then it is not hard to check that, for any fixed $A\in e_7^{\C}$,
the following map
\begin{eqnarray}
\Lambda^2 V &\lon & \rho(e_7^{\C} + \C) \nonumber\\
u\ot v        & \lon & 2\lambda \mu <\!u,v\!> \rho(A)
+ \rho(u\circ \rho(A)v) -  \rho(v\circ \rho(A)u)  
\label{K}
\end{eqnarray}
defines an element of $\rho(e_7^{\C} + \C)\ot \Lambda^2 V^*$ which lies
in the kernel of the composition
$$
\rho(e_7^{\C} + \C)\ot \Lambda^2 V^*\lon V\ot V^* \ot \Lambda^2 V^*
 \lon V\ot \Lambda^3 V^*.
$$
Thus, the above formula gives an explicit realization of the
isomorphism $K(\rho(e_7^{\C} + \C))=\Ad(e_7^{\C})$. In particular,
it shows that $K(\rho(e_7^{\C} + \C))=K(\rho(e_7^{\C}))$.

Having obtained an explicit structure of $K(\rho(e_7^{\C}))$, it is
straightforward to show that a generic element of 
$K^1(\rho(e_7^{\C})) \subset \rho(e_7^{\C})\ot V^*\ot \Lambda^2 V^*$
is of the form
\begin{eqnarray*}
V\ot \Lambda^2 V &\lon & \rho(e_7^{\C}) \nonumber\\
s\ot u\ot v        & \lon & 2\lambda\mu <\!u,v\!> \rho(s\circ w) + 
\rho(u\circ \rho(s\circ w)v) -  \rho(v\circ \rho(s\circ w)u)  
\end{eqnarray*}
for some fixed $w\in V\simeq V^*$. This establishes the isomorphism
$K^1(\rho(e_7^{\C}))=V^*$.  $\Box$

\begin{center}
{\bf \S 3 A construction of torsion-free connections}
\end{center}

We briefly describe here the construction of torsion-free connections
with prescribed holonomy which
was presented in \cite{CMS}.

Let $\g \subset \gl(V)$ a Lie sub-algebra where $V$ is a finite-dimensional
vector space. 

A $G$-equivariant $C^\infty$-map $\phi: \g^* \rar \Lambda^2 V^*$ is called
{\em admissible} if for every $p \in \g^*$, the map $d\phi_p^*: \Lambda^2 V
\lon \T_p^*\g^* \simeq \g$ lies in $K(\g)$.

For a given admissible map $\phi$, one may define the following Poisson
structure on the dual $W^*$ of the semi-direct Lie algebra $W = \g \oplus V$:
\begin{equation} \label{Poisson}
\{f,g\}_{p+\nu} = p ([A,B]) + \nu(A \cdot y - B \cdot x) + \phi(p)(x,y),
\end{equation}
where $df = A + x$ and $dg = B + y$ are the decompositions of
$df, dg \in T^*W^* \simeq \g \oplus V$, and $p \in \g^*, \nu \in V^*$. This
Poisson structure
may be regarded as a deformation of the natural linear Poisson structure on
$W^*$.

Let $\pi: S \rar U$ be a {\em symplectic realization} of an
open subset $U \subset W^*$, i.e. $\pi$ is a submersion from a symplectic
manifold $S$ with 
symplectic 2-form $\Omega$ 
such that
$$
\{\pi^*(f), \pi^*(g)\}_S = \pi^*(\{f,g\})  \ \ \  \mbox{for all $f,g \in
C^\infty(P,\R)$,}
$$
where $\{\ ,\ \}_S$ is the Poisson structure on $S$ induced by the
symplectic structure. At those points where the rank of the Poisson
structure is maximal, such a symplectic realization exists at least locally.

Regarding each element $w \in W \simeq T^* W^*$ as a 1-form
on $W$, we define the distribution 
$$
\cD= \{ \xi_w := \#(\pi^*(w)) \mid w\in W\} \subset TS
$$
on $S$, where $\#$ is the index-raising map induced by $\Omega$.
Since $\Omega$ is non-degenerate, \linebreak $\rank \cD = \dim W$. Moreover,
for the bracket relations one calculates

\begin{equation}\label{brackets}
\begin{array}{llll} [\xi_A, \xi_B] = \xi_{[A,B]}, & [\xi_A, \xi_x] =
\xi_{A \cdot x}, & \mbox{and} & [\xi_x, \xi_y](s) = \xi_{d\phi_p^* (x,y)},
\end{array}
\end{equation}
where $A,B \in \g$, $x,y \in V$ and $p = \pi(s)$.

Let $F\subset S$ be an integral leaf 
of $\cD$. By the very definition of $\cD$, $F$ comes equipped with a
$W$-valued coframe $\theta + \om$, where $\theta$ and $\om$ take values
in $V$ and $\g$ respectively, defined by the equation
$$
(\om + \theta) (\xi_w) = w.
$$

Note that by the first equation in (\ref{brackets}), the vector fields
$\xi_A$, $A \in \g$, induce a free local group action of $G$ on $F$, where
$G \subset Gl(V)$ is the connected Lie subgroup corresponding to $\g \subset
\gl(V)$. After shrinking $F$ as necessary, we may assume that $M := F/G$ is
a manifold. Standard arguments then imply that there is a unique embedding
of $\imath: F \hookrightarrow {\frak F}_V$, where ${\frak F}_V$ denotes the
$V$-valued coframe bundle of $M$ and a torsion-free connection on $M$ such
that $\imath^*(\underline\theta + \underline\om) = \theta + \om$, where
$\underline\theta$ and $\underline\om$ denote the tautological and the
connection 1-form on ${\frak F}_V$, respectively. Clearly, the holonomy of
this connection is contained in $G$; in fact, by the {\em
Ambrose-Singer-Holonomy Theorem}, the holonomy algebra is generated by
$\{ d\phi^*_p(x,y)\ |\ x,y \in V, p \in \pi(F)\}$. A connection which comes
from this construction is called a {\em Poisson connection}. This leads to
the following

\begin{theorem}{\em \cite{CMS}}\label{T2}
Let $\g \subset \gl(V)$ be a
Lie sub-algebra where $V$ is a finite-dimensional vector space, and
let
$$
K_0(\g) = \left\{ R \in K(\g) \mid \mbox{\em span}\{R(x,y), \ \mbox{\em all}\
x,y\in V\}= \g\right\}.
$$
If $\phi: \g^* \rar
\Lambda^2 V^*$ is admissible, and if the open set $U_0 \subset
\g^*$ given by
\[
U_0 := (d\phi^*)^{-1} (K_0(\g))
\]
is non-empty, then there exist Poisson connections induced by $\phi$
whose holonomy representations are equivalent to $\g$. Moreover, if
$\phi|_{U_0}$ is not affine, then not all of these connections are
locally symmetric.
\end{theorem}

It is not clear at present how general the class
of Poisson connections is, nor how many irreducible
Lie subalgebras $\g\subset \gl(V)$ admit admissible
maps $\phi: \g^*\rar \Lambda^2V^*$ which are not affine. However, there is a
class of Lie subalgebras for which the above construction exhausts {\em all
possible torsion-free connections} with this holonomy. Namely, we define
the $\g$-module 
$$
{\cal P}^{(1)}(\g) = (\odot^2 \g\ot \Lambda^2 V^*) \cap (\g\ot K(\g))
\subset V\ot K^1(\g)
$$
and regard elements $\phi_2 \in {\cal P}^{(1)}(\g)$ as polynomial maps
$\g^* \rar \Lambda^2 V^*$ of degree $2$. It is then obvious that each
$G$-invariant $\phi_2 \in {\cal P}^{(1)}(\g)$ is admissible, and we have
the following result.

\begin{theorem}{\em \cite{CMS}}
\label{T3} Let $\g \subset \gl(V)$ be an irreducibly acting subalgebra, and
suppose that there is an invariant element 
$\phi_2 \in {\cal P}^{(1)}(\g)$
such that the associated $G$-equivariant linear maps
\begin{equation} \label{isomorphisms}
\begin{array}{rrcl}
\phi_2': & \g^* & \longrightarrow & K(\g)\\
\phi_2'': &  V^* & \longrightarrow & K^1(\g).
\end{array} 
\end{equation}
are isomorphisms. Then every torsion-free affine connection whose holonomy
algebra is contained in $\g$ is a Poisson connection
induced by an admissible map
$$
\phi = \phi_2 + \tau,
$$
where $\tau \in \Lambda^2 V^*$ is a (possibly vanishing)
$\g$-invariant 2-form. In particluar, the moduli space
of such connections is finite dimensional, and each such connection is
analytic. Also, the dimension of the symmetry group of this connection
equals $\dim W^* - 2k$ where $k$ is the half-rank of the Poisson structure
on $W^*$ induced by $\phi$ in (\ref{Poisson}).
\end{theorem}

At first sight, the premise that the maps (\ref{isomorphisms}) be
isomorphisms looks like an unreasonably strong condition in order to
utilize this Theorem. Nevertheless, this premise {\em does} hold for the
exotic holonomies $\SO(p,q)\SL(2,\R)$ and $\SO(n,\C)\SL(2,\C)$ which were
discovered in \cite{CMS}. Also, we will show in \S 4 that it
also holds for the representations $E_7^{(a)}$ from the main
theorem.

For the proof, we shall need the following version of Schur's Lemma:

\begin{lemma} \label{L4} Let $\g$ be a reductive Lie algebra, and suppose
that $\g$ acts irreducibly on the finite dimensional vector spaces $V$ and
$W$. If $\rho: V \rar W$ is a linear map satisfying \begin{eqnarray}
\label{Schur} A\ \rho\ B = B\ \rho\ A & \mbox{for all $A, B \in \g$,}
\end{eqnarray} then $\rho = 0$.
\end{lemma}

\medskip \noindent {\bf Proof of Theorem \ref{T3}} Let $F \subset 
{\frak F}_V$ be a
$G$-structure on the manifold $M$ where ${\frak F}_V \rar M$ is the
$V$-valued coframe bundle of $M$, and denote the tautological $V$-valued
1-form on $F$ by $\theta$. Suppose that $F$ is equipped with a torsion-free
connection, i.e. a $\g$-valued 1-form $\om$ on $F$. Since $\phi_2'$ is
an isomorphism, the {\em first and second structure equations} read
\begin{equation} \label{eq:struct} \begin{array}{ll}
d\theta & = - \om \wedge \theta\\
d\om & = - \om \wedge \om - 2 (\phi_2'({\bf a})) \circ (\theta \wedge
\theta),\end{array} \end{equation}
where ${\bf a}: F \rar \g^*$ is a $G$-equivariant map. Differentiating
(\ref{eq:struct}) and using that $\phi_2''$ is an isomorphism yields
the {\em third structure equation} for the differential of ${\bf a}$:
\begin{equation} \label{eq:struct3}
d{\bf a} = -\om \cdot {\bf a} + \jmath({\bf b} \otimes \theta),
\end{equation}
for some $G$-equivariant map ${\bf b}: F \rar V^*$, where $\jmath:  V^*
\otimes V \rar \g^*$ is the natural projection. The multiplication in
the first term refers to the coadjoint action of $\g$ on $\g^*$.

Let us define the map ${\bf c}: F \rar V^* \otimes V^*$ by
\begin{equation} \label{eq:struct4}
{\bf c}_p (x,y) := d{\bf b}(\xi_x) (y) - \phi_2({\bf a}_p, {\bf a}_p, x, y).
\end{equation}
Differentiation of (\ref{eq:struct3}) yields
\begin{equation} \label{eq:cinvariance} \begin{array}{ll}
{\bf c}_p(x, Ay) = {\bf c}_p(y, Ax) & \mbox{for all $x,y \in V$ and all $A
\in \g$.}
\end{array} \end{equation}

If we let $\rho: V \rar V^*, x \mapsto {\bf c}_p(x,\underline{\ \ }) +
{\bf c}_p(\underline{\ \ }, x)$, then (\ref{eq:cinvariance}) and Lemma
\ref{L4} imply that $\rho = 0$, i.e. ${\bf c}_p$ must be skew-symmetric and
$G$-invariant. This and differentiation of (\ref{eq:struct4}) implies that
\[
d{\bf c} = 0,
\]
i.e. ${\bf c}_p \equiv \tau \in \Lambda^2 V^*$ is {\em constant}. Thus, the
$G$-equivariance of ${\bf b}$ and (\ref{eq:struct4}) yield
\begin{equation} \label{eq:newstruct4}
d{\bf b} = -\om \cdot {\bf b} + \left( {\bf a}_p^2 \contr \phi_2 + \tau
\right) \circ \theta,
\end{equation}
where $\contr$ refers to the contraction of ${\bf a}_p^2 \in \odot^2(\g^*)$
with $\phi_2 \in \odot^2(\g) \otimes \Lambda^2 V^*$.

Let us now define the Poisson structure on $W^* = \g^* \oplus V^*$
induced by $\phi := \phi_2 + \tau$, and let $\pi := {\bf a} + {\bf b}: F \rar
W^*$. From (\ref{eq:struct3}) and (\ref{eq:newstruct4}), one can now show
that, at least locally, the connection is indeed a Poisson connection
induced by $\phi$.

Let ${\frak s} \subset {\frak X}(F)$ be the Lie algebra of infinitesimal
symmetries. Let $f: W^* \supset U \rar {\Bbb F}$ be a local function
which is
constant on the symplectic leaves. Then it is easy to see that $\#
\pi^*(df)$ is an infinitesimal symmetry. It follows that $\dim {\frak s}
\geq \dim W^* - 2k$. On the other hand, if $X \in {\frak s}$ then
$\pi_*(X) = 0$, hence $\dim {\frak s} \leq \dim W^* - 2k$.

The statements about analyticity and the moduli space are now immediate.
{\hfill \rule{.5em}{1em}\mbox{}\bigskip}

\medskip \noindent {\bf Proof of Lemma \ref{L4}} Throughout the proof, we
make the simplifying assumption that $\rank \g >1$, as the case $\rank \g
= 1$ is straightforward. Let $P \subset V^* \ot W$ be the subspace of all
maps $\rho: V \rar W$ satisfying (\ref{Schur}). It is easy to verify that
$P$ is $\g$-invariant. We complexify $\g$, $V$ and $W$ and pick Cartan and
weight space decompositions \[ \begin{array}{llll} \g = {\frak t} \oplus
\bigoplus_\alpha \g_\alpha, & V = \bigoplus_{\mu} V_\mu &
\mbox{and} & W = \bigoplus_{\mu} W_\mu. \end{array} \]
Let $\rho \in P$, and let $x_\mu \in V_\mu$ with $\mu \neq 0$. Then choosing
$A, B \in {\frak t}$, with $\mu(A) = 0,\ \mu(B) \neq 0$, (\ref{Schur})
implies that $A \rho x_\mu = 0$, and therefore,

\begin{equation} \label{**} \rho x_\mu \in \sum_k W_{k \mu},
\end{equation}
where the sum is taken over all weights of $W$ which are scalar multiples
of $\mu$. Now let $\rho_\lambda \in P$ be an element of weight
$\lambda \neq 0$. Then
$\rho_\lambda x_\mu \in W_{\lambda + \mu}$, and thus from (\ref{**}) we
conclude:
\begin{equation} \label{indep}
\mbox{$\rho_\lambda x_\mu = 0$ whenever $\lambda, \mu$ are linearly
independent.}
\end{equation}

Let $x_{k \lambda} \in V_{k \lambda}$ be a weight vector with $k \neq 0$,
and let $A_\alpha \in \g_\alpha$ where $\alpha$ is a root independent of
$\lambda$. Then, using (\ref{indep}) twice, we get

\[ \begin{array}{ll}
0 & = (A_\alpha \cdot \rho_\lambda) x_{k \lambda}\\
  & = A_\alpha (\rho_\lambda x_{k \lambda}) - \rho_\lambda
      (A_\alpha x_{k \lambda})\\
  & = A_\alpha (\rho_\lambda x_{k \lambda}) - 0.
\end{array}
\]

Next, note that $V_0$ is spanned by elements of the form $A_\alpha
x_{-\alpha}$ with $A_\alpha \in \g_\alpha$ and $x_{-\alpha} \in
V_{-\alpha}$. If $\alpha, \lambda$ are independent, pick $A_0 \in
{\frak t}$ with $\alpha(A_0) = 0$ and $\lambda(A_0) \neq 0$. Then
(\ref{Schur}) implies that $\rho_\lambda (A_\alpha x_{-\alpha}) = 0$.
Finally, if $\alpha, \lambda$ are dependent, then for $\beta \neq \pm
\alpha$, we get from (\ref{Schur}) and (\ref{indep}) that $A_\beta
(\rho_\lambda (A_\alpha x_{-\alpha})) = 0$.

In either case, we get that for any $\mu$, $A_\alpha (\rho_\lambda x_\mu)
= 0$ whenever $\lambda, \alpha$ are independent, and hence

\[ \g \cdot (\rho_\lambda V) \subset \sum_k W_{k \lambda}. \]
Since there must be weights in $W$
independent of $\lambda$, and since $W$ is irreducible, we conclude that
$\rho_\lambda = 0$, contradicting $\lambda \neq 0$.

Thus, $P$ has no weights $\neq 0$, i.e. $P$ is acted on trivially by $\g$,
and from there it is easy to conclude that $P = 0$. $\Box$

\begin{center}
{\bf \S 4 Proof of the main theorem}
\end{center}

Let $\g \subset \gl(V)$ be one of the representations in the Main Theorem.
Evidently, the Main Theorem will follow from Theorems~\ref{T2} and \ref{T3}
if we can find an element $\phi_2 \in {\cal P}^{(1)}(\g)$ such that
$K_0(\g)$ is dense in
$K(\g)$, and the corresponding maps in (\ref{isomorphisms}) are
isomorphisms. In particular, (iv) of the main Theorem follows since, in
each case, $\dim W^* = \dim V + \dim \g = 56 + 133$ is odd.

The density of $K_0(\g)$ in $K(\g)$ follows immediately from (\ref{K}).

To compute ${\cal P}^{(1)}(\g)$, first note that the $\g$-module 
$K(\g)\ot \g\simeq \g\ot \g$ has only one 1-dimensional
$\g$-submodule. So, if there is an invariant element $\phi_2$ in
${\cal P}^{(1)}(\g)$,
it is unique up to a non-zero scalar factor. Since $\phi_2': \g^* \lon
K(\g)$ must be an isomorphism, the formula (\ref{K}) leaves no choice
but the following element of  $\g\ot K(\g)$ as a candidate
for $\phi_2$:
$$
\phi_2(C,D, u,v) = 2\lambda\mu<\!u,v\!>B(C,D) + B(u\circ \rho(C)v, D)  -
B(v\circ \rho(C)u, D), 
$$
where $C,D\in \g$, $u,v\in V$ and where we identify $\g=\g^*$ via the
Killing form. Clearly, this element is $\g$-invariant. 
Since $B(u\circ \rho(C)v, D) = \lambda <\!\rho(D)u, \rho(C)v\!>$,
we have 
\begin{eqnarray*}
\phi_2(C,D, u,v) & = & 2\lambda\mu<\!u,v\!>B(C,D) + \lambda <\!\rho(D)u, 
\rho(C)v\!> - \lambda <\!\rho(D)v, \rho(C)u\!>  \\
&=&  2\lambda\mu<\!u,v\!>B(C,D) - \lambda <\!\rho(C)v, \rho(D)u\!>
- \lambda <\!\rho(D)v, \rho(C)u\!>  \\
\end{eqnarray*}
which makes it evident that $\phi_2\in (\odot^2\g\ot \Lambda^2 V^*)
 \cap (\g\ot K(\g))$.
That $\phi_2': \g^* \lon K(\g)$ is an isomorphism follows from the
very  definition of $\phi_2$. Since $\phi_2'': V^* \lon K^1(\g)$ is
evidently non-zero, then, by Theorem~\ref{T1} and the $\g$-invariance
of $\phi_2$, $\phi''_2$ must be an isomorphism as well.
$\Box$
\vspace{5 mm}

\noindent{\small {\em Acknowledgement}. One of the authors (SM)
would like to thank Andrew Swann for valuable discussions and
for drawing his attention to the very informative lecture notes
\cite{Adams}.}

\pagebreak

\[ \small \begin{array}{ll}
\mbox{\sc Quo-Shin Chi} &
\mbox{\sc Department of Mathematics, Campus Box 1146,}\\ &
\mbox{\sc Washington University, St. Louis, Mo 63110, USA}\\ &
\mbox{\rm chi@artsci.wustl.edu} \\

\mbox{\sc Sergey A.  Merkulov} &
\mbox{\sc Department of Pure Mathematics, Glasgow Univer-}\\ &
\mbox{\sc sity, 15 University Gardens, Glasgow G12 8QW, UK}\\ &
\mbox{\rm sm@maths.glasgow.ac.uk} \\

\mbox{\sc Lorenz J.  Schwachh\"{o}fer} &
\mbox{\sc Mathematisches Institut, Universit\"at Leipzig,}\\ &
\mbox{\sc Augustusplatz 10-11, 04109 Leipzig, Germany}\\ &
\mbox{\rm schwachhoefer@mathematik.uni-leipzig.de}
\end{array}
\]


{\small
\begin{thebibliography}{99}

\bibitem{Adams} J.F.\ Adams, {\em Exceptional Lie groups}, unpublished
lecture notes.

\bibitem{Ah} D.N, Ahiezer, D.N., {\em   Homogeneous complex
manifolds}, in: Several Complex Variables IV, Springer 1990.

\bibitem{BE} {R.J.\ Baston and M.G.\ Eastwood}, {\em The Penrose 
transform, its interaction with representation theory}, Oxford University 
Press, 1989.

\bibitem{Berger} { M.\ Berger}, {\em Sur les groupes d'holonomie des
vari\'{e}t\'{e}s \'{a} connexion affine et des vari\'{e}t\'{e}s
Riemanniennes},  {Bull.\ Soc.\ Math.\ France}\, {\bf 83} (1955), 279-330.


\bibitem{Bryant2} {R.\ Bryant}, {\em Two exotic holonomies in 
dimension four, path geometries, and twistor theory}, {Proc.\ Symposia in
Pure Mathematics}\, {\bf 83} (1991), 33-88.

\bibitem{Bryant3} {R.\ Bryant}, {\em Classical, exceptional,
and exotic holonomies: a status report}, preprint 1995.

\bibitem{CS} {Q.-S.\ Chi, L.J.\ Schwachh{\"o}fer}, {\em Exotic holonomy
on moduli spaces of rational curves}, preprint

\bibitem{CMS} {Q.-S.\ Chi, S.A.\ Merkulov, L.J.\ Schwachh\"{o}fer},
{\em On the existence of infinite series of exotic holonomies},
to appear in Inv.\ Math.

\bibitem{HS}  {M.\ Hausner and J.T.\ Schwartz}, {\em Lie groups,
Lie algebras}, Gordon and Breach 1968

\bibitem{Me1} S.A.\ Merkulov,  Addendum in \cite{Manin}.


\bibitem{Manin} Yu.I.\ Manin, {\em Gauge field theory and complex
geometry}, 2nd Edition, Springer-Verlag, in press.

\bibitem{On} A.L.\ Onishchik, E.B.\ Vinberg, {\em Lie groups and
algebraic groups}, Springer-Verlag 1990


\end{thebibliography}
}
\end{document}